# Optimal perturbations for nonlinear systems using graph-based optimal transport

Piyush Grover[1]

*Mitsubishi Electric Research Labs, Cambridge, MA, USA*

Karthik Elamvazhuthi

*Arizona State University, Tempe, AZ, USA*

**Abstract**

We formulate and solve a class of finite-time transport and mixing problems in the set-oriented framework. The aim is to obtain optimal discrete-time perturbations in nonlinear dynamical systems to transport a specified initial measure on the phase space to a final measure in finite time. The measure is propagated under system dynamics in between the perturbations via the associated transfer operator. Each perturbation is described by a deterministic map in the measure space that implements a version of Monge-Kantorovich optimal transport with quadratic cost. Hence, the optimal solution minimizes a sum of quadratic costs on phase space transport due to the perturbations applied at specified times. The action of the transport map is approximated by a continuous pseudo-time flow on a graph, resulting in a tractable convex optimization problem. This problem is solved via state-of-the-art solvers to global optimality. We apply this algorithm to a problem of transport between measures supported on two disjoint almost-invariant sets in a chaotic fluid system, and to a finite-time optimal mixing problem by choosing the final measure to be uniform. In both cases, the optimal perturbations are found to exploit the phase space structures, such as lobe dynamics, leading to efficient global transport. As the time-horizon of the problem is increased, the optimal perturbations become increasingly localized. Hence, by combining the transfer operator approach with ideas from the theory of optimal mass transportation, we obtain a discrete-time graph-based algorithm for optimal transport and mixing in nonlinear systems.

## 1. Introduction

In the study of nonlinear dynamical systems, the problem of efficient phase space transport is of central interest. For extraction of organizing phase-space structures in low-dimensional systems arising in fluid kinematics [1, 2], celestial mechanics [3], and plasma physics [4], several methods based on geometric [5], topological [6, 7] and statistical techniques [8, 9, 10, 11] have been developed.

The geometric techniques focus on extracting the Lagrangian coherent structures in autonomous and non-autonomous systems, which are often the stable and unstable manifolds [12] of fixed points or periodic orbits, or their time-dependent analogues [13, 14]. Techniques based on lobe-dynamics [12] allow for quantifying transport between different weakly mixing regions in the phase space ('coherent sets'). Once these structures have been identified, intelligent control strategies can be formulated to obtain efficient phase space transport between the coherent sets in the phase space [3, 15, 16]. Based on these geometric methods, a method for computing efficient chaotic transport was described in Ref. [17]. Furthermore, methods for computing impulsive perturbations to differential equations for purpose of optimal enhancement of mixing have also been developed recently [17, 18].

---

[1]Corresponding Author. Email: grover@merl.com

Based on a class of mixing measures which capture the large scale inhomogeneities in the scalar field [19], locally optimal switching laws among a finite set of optimal velocity fields were obtained via optimal control techniques in Ref. [20]. Via analytic computation of derivatives of mixing norms, local-in-time optimal fields among an energy or enstrophy constrained set of incompressible velocity fields were obtained in Ref. [21]. This work inspired a series of works on obtaining bounds on mixing rates [22, 23, 24] based on various constraints on the advection fields.

Statistical set-oriented methods for computing transfer operators [9, 25] allow the discovery of 'coherent sets' in autonomous and non-autonomous dynamical systems. Furthermore, rigorous optimal control methods using set-oriented methods have also been developed [26, 27]. In Refs. [28, 29, 30], an optimal control framework for asymptotic stabilization of arbitrary initial measure to an attractor is presented. This framework is based on computing a (control) 'Lyapunov measure', which is a measure-theoretic analogue of control Lyapunov function. Also relevant is the work in the area of occupation measures, see Ref. [31].

In this paper, we are interested in the problem of discrete-time 'optimal transport' [32] under nonlinear dynamics. Given an initial measure in phase space, the problem can be stated as follows : *Compute the sequence of optimal perturbations applied to the measure at discrete times, such that a desired measure at prescribed final time is achieved.* Here, the optimality implies minimizing an appropriate cost on the phase space transport that occurs due to the perturbations. Between the discrete times, the measure evolves under the action of given nonlinear system dynamics. If the desired final measure is chosen to be uniform over a compact phase space, the problem becomes that of *optimal enhancement of mixing*. These problems are meaningful if the perturbations are 'small' compared to the underlying dynamics, in some appropriate norm. These problems are motivated by several potential applications. The first application, similar to several studies mentioned earlier, is to transport and mixing of scalars in fluid systems. In this case, if the control can be applied at a faster time-scale than the underlying dynamics, modeling the perturbations as instantaneous is a good approximation. Another motivation comes from the problem of controlling swarms of agents with similar dynamics in an ambient flow field. For instance, the control of magnetic particles in blood stream [33, 34, 35], robotic miniature bees in air [36, 37], and swarms of autonomous underwater vehicles (AUVs) in the ocean [38] can all be studied as swarm control and planning problems in presence of an ambient flow field. Here, one represents the distribution of swarms in phase space by measures, and the control problem can be formulated in terms of measure transport.

In Ref. [39], the problem of computing local perturbations for enhancing mixing is addressed using statistical methods. A set-oriented transfer operator approach is used to propagate the dynamics, and perturbations are modeled via a stochastic kernel. The resulting convex optimization problem is then solved for perturbations that lead to minimum difference in $L^2$ norm from the desired density at each time step. In Ref. [40], an infinitesimal generator based approach is employed to increase rate of mixing in the continuous-time setting. Hence, the previous work has either focussed on maximizing mixing rates under energy, enstrophy or other constraints, or maximizing the 'mixedness' of *final* state under constraints.

We study a different class of finite-time transport (and mixing) optimization problems, where the aim is to minimize the sum of phase space transport cost due to perturbations applied at specified times. We retain the set-oriented transfer operator framework, and the use of discrete-time perturbations, introduced in context of mixing optimization in Ref. [39], and adapt it to the problem of interest. First, we use a set-oriented version of an optimal transport cost, motivated by the theory of optimal mass transportation [32], as the cost of the perturbations. This cost can be understood as a quadratic cost of phase space transport due to perturbations, and it also has a control theoretic interpretation [41]. Second, our problem is formulated as fixed final time problem with prescribed initial and final measures. The maps corresponding to discrete-time perturbations are computed via continuous pseudo-time advection between intermediate measures supported on a graph. This graph is built upon the same box discretization used for computing the transfer operator. Hence, we formulate and solve a convex global-in-time graph-based optimization problem that switches between flow due to the dynamics of the system, and a pseudo-time advection due to perturbations, while minimizing the optimal transport cost. Optimal mass transportation is concerned with optimization of measure transport under different settings, and has deep connections with phase space transport in dynamical systems [42, 43, 44]. Hence, it forms a natural setting in which to study the problem of interest. Since set-oriented methods have been very successful in the study of complex dynamical systems, it is desirable that the optimization framework can exploit the computational advantages provided by such methods. The graph-based optimal mass transportation method that we use is a natural tool for the set-oriented framework.

This paper is organized as follows. In Section 2, we give a brief overview of the theory of optimal mass transport, and a continuous time approach to solving the canonical optimal transport problem. In Section 3, we formally state the optimization problem, and describe our computational framework in detail. In Section 4, we apply our method to problems of transport and mixing in a well-studied fluid dynamical system, the time-periodic double-gyre. Finally, in Section 5, we provide conclusions and discuss some avenues for future research.

## 2. Optimal Transport and Dynamical Systems

### 2.1. Monge-Kantorovich problem

The Monge-Kantorovich optimal transport (OT) problem [32] is concerned with transport of an initial measure $\mu_0$ on a space $X$ to a final measure $\mu_1$ on a space $Y$. In the original formulation, it involves solving for a measurable transport map $T : X \to Y$, which pushes forward $\mu_0$ to $\mu_1$ in an optimal manner. The cost of transport per unit mass is prescribed by a function $c(x, T(x)) : X \times M(X, Y) \to [0, \infty)$, where $M(X, Y)$ is the space of all measurable maps from $X$ to $Y$. Hence, the optimization problem is to compute

$$\min_T \int c(x, T(x))d\mu_0(x), \tag{1}$$
$$\text{s.t. } T_\# \mu_0 = \mu_1,$$

where $T_\#$ is the pushforward of $T$, i.e. $(T_\# \mu)(A) = \mu(T^{-1}(A))$ for every $A$. The square root of optimal cost obtained as solution of Eq. (1) with $c(x, y) = \|x - y\|^2$ is called the 2−Wasserstein distance, and we denote it by $W_2(\mu_0, \mu_1)$.

### 2.2. Benamou-Brenier fluid dynamics approach

The OT problem Eq. (1) described in the previous section is concerned with only the measures at initial and final time. For $c(x, y) = \|x - y\|^2$, the optimal map $T$ is known to be of the form $T(x) = \nabla \psi(x)$ for some convex function $\psi(x)$ (see Section 1.3 in Ref. [45]). Then, the pushforward constraint can written as,

$$det(D^2 \psi(x))\mu_1(\nabla \psi(x)) = \mu_0(x) \tag{2}$$

Numerically solving the optimization problem with nonlinear constraint in Eq. (2) is difficult. An alternative approach, inspired by fluid dynamics, was described in Ref. [46]. Under some additional conditions on smoothness of measures [45], the optimization problem is formulated in terms of an advection field $u(x, t)$, and initial and final *densities* $(\rho_0(x), \rho_1(x))$. The core idea is to obtain the optimal map $T$ as a result of advection over a time period $(t_0, t_f)$ by the optimal advection field $u(x, t)$. It can be shown that the optimization problem given by Eq. (1) (with $X = Y = \mathbb{R}^\mathbf{d}$) for the case of quadratic cost, $c(x, y) = \|x - y\|^2$, is equivalent to the following problem:

$$W_2(\mu_0, \mu_1)^2 = \min_{u(x,t)} t_f \int_{\mathbb{R}^n} \int_{t_0}^{t_f} \rho(x,t)|u(x,t)|^2 dt dx, \tag{3}$$
$$\text{s.t. } \frac{\partial \rho(x,t)}{\partial t} + \nabla \cdot (\rho(x,t)u(x,t)) = 0 \ \ t_0 \leq t \leq t_f,$$
$$\rho(x, t_0) = \rho_0(x), \rho(x, t_f) = \rho_1(x).$$

This aim of this problem can be understood as minimizing the time integral of the total kinetic energy of the 'fluid' being advected by the field $u(x, t)$, subjected to initial and final densities. Furthermore, the optimal advection field is a potential flow, i.e., $u(x,t) = \nabla \phi(x,t)$ for some potential field $\phi(x,t)$ [46]. The Euler-Lagrange equations for this optimization problem can be written as

$$\frac{\partial \phi}{\partial t} + \frac{|\nabla \phi|^2}{2} = 0, \tag{4}$$

which are pressure-less version of Euler equations.

By a change of variables from $(\rho, u)$ to $(\rho, m \stackrel{\Delta}{=} \rho u)$, the optimization problem in Eq.(3) can be put into a convex form. The transformed optimization problem is

$$\min_{m(x,t)} t_f \int_{\mathbb{R}^n} \int_{t_0}^{t_f} \frac{|m(x,t)|^2}{\rho(x,t)} dt dx, \tag{5}$$

$$\text{s.t. } \frac{\partial \rho(x,t)}{\partial t} + \nabla \cdot (m(x,t)) = 0 \ t_0 \leq t \leq t_f,$$

$$\rho(x, t_0) = \rho_0(x), \rho(x, t_f) = \rho_1(x).$$

In the above equations, the constraints are now linear in the problem variables $(\rho, m)$. The term inside the integral in the cost function, $\frac{|m(x,t)|^2}{\rho(x,t)}$, is of the 'quadratic-over-linear' form [47], and can be shown to be convex [45] w.r.t $(\rho, m)$. Hence, by transforming the transport problem into continuous time, and by using a change of variables, the non-convex problem in Eq. (1) for $c(x,y) = \|x - y\|^2$, has been converted into a convex problem in Eq. (5).

## 3. Problem Setup and Computational Approach

Consider a map $F : X \to X$ on phase space $X$. This map may be obtained from a discrete time dynamical system, or as a time-$t$ map of the flow of an autonomous or time-periodic dynamical system. Given a pair of initial and final desired measures $(\mu_{t_0}, \mu_{t_f})$ on $X$, we want to obtain optimal measure-preserving discrete-time perturbations acting on the measure, represented by linear operators $\tilde{T}^1, \tilde{T}^2, \ldots, \tilde{T}^n$ applied at $n$ discrete time-steps. The $\tilde{T}^i$s can be thought of as pushforwards of perturbations that act on the phase space. The measure evolves under the action of $F$ in between the time steps. The cost function, evaluated between the pair of measures before and after the application of each $\tilde{T}^i$, should reflect a quadratic cost of phase transport due to perturbations, in the spirit of discussion in Sec. 2 . We motivate and develop an appropriate formulation for this cost in this section. This is a fixed final-time problem with fixed final 'state', namely the measure $\mu_{t_f}$. The problem can be formally stated as,

$$\min_{T^1, T^2, \ldots, T^n} \sum \|\mu^i - \hat{\mu}^i\|, \tag{6}$$

$$\hat{\mu}^i = F_\# \mu^{i-1} \quad \forall i \in (1, 2, \ldots, m), \tag{7}$$

$$\mu^i = \tilde{T}^i \hat{\mu}^i \quad \forall i \in (1, 2, \ldots, m), \tag{8}$$

$$\mu^0 = \mu_{t_0}, \mu^n = \mu_{t_f}, \tag{9}$$

where the action of each $\tilde{T}^i$ is measure preserving. In the measure space, the uncontrolled dynamics of a nonlinear system are described by a linear operator, called the 'transfer operator' or 'Perron-Frobenius operator' [10, 48]. We model the discrete-time perturbations as deterministic maps in this measure space, obtained by solving a continuous (pseudo-) time Monge-Kantorovich problem. We obtain a convex optimization problem by switching between the uncontrolled transport and perturbations. For computation in the infinite dimensional space, we take the 'discretize-then-optimize' approach.

## 3.1. Transfer Operator

The Perron-Frobenius transfer operator $P$ [10] corresponding to the map $F$, is a linear operator which pushes forward measures in phase space according to the dynamics of the trajectories under $F$. Let $\mathbf{B}(X)$ denote $\sigma$−algebra of Borel sets in $X$. Then,

$$P\mu(A) = \mu(F^{-1}(A)) \quad \forall A \in \mathbf{B}(X). \tag{10}$$

Here $F$ is assumed to be non-singular w.r.t. the Lebesgue measure. The transfer operator lifts the evolution of the dynamical systems from phase space $X$ to the space of measures $\mathbf{M}(X)$. Numerical approximation of $P$, denoted by $\hat{P}$, may be viewed as a transition matrix of an $N$-state Markov chain [11]. For computation, we partition the phase space volume of interest into $N$ boxes, $B_1, B_2, \ldots, B_N$. $\hat{P}$ is computed via the Ulam-Galerkin method [49, 11], as follows

$$\hat{p}_{ij} = \frac{\bar{m}\left(B_i \cap F^{-1}(B_j)\right)}{\bar{m}(B_i)}, \tag{11}$$

where $\bar{m}$ is the Lebesgue measure. The process is illustrated in Fig. 1.

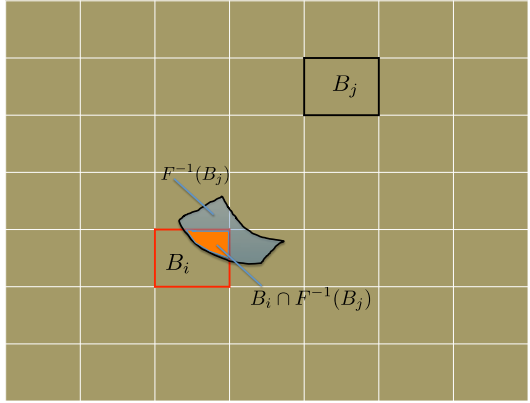

**Figure 1:** Computation of the transition matrix $\hat{P}$ by a set-oriented method. Box $B_j$ at the current time instance is mapped (backwards) to $F^{-1}(B_j)$ at the previous time instance. The value of the entry $\hat{p}_{ij}$ is the fraction of box $B_i$ that is mapped into box $B_j$ by $F$.

## 3.2. Monge-Kantorovich Transport on graphs

Now consider a graph $G = (V, E)$ on $X$, where the set of vertices $V$ represent the boxes $B_i$, and the set of directed edges $E$ are obtained from the topology of $X$. Let us consider formulating a quadratic cost optimal transportation problem for measures $(\mu_0, \mu_1)$ supported on $V$, similar to the original Monge-Kantorovich formulation. The problem will involve computing a measurable transport map $\tilde{T}(v, w)$ between each pair of vertices $(v, w)$, such that some quadratic cost is minimized. If $\bar{d}(v, w)$ denotes the shortest path distance between the two vertices, then a reasonable candidate for quadratic transport cost in Eq. (6), analogous to the one given in Eq. (1) with $c(x, y) = \|x - y\|^2$, is

$$W_{2,N}(\mu_0, \mu_1)^2 = \min_{\tilde{T}} \sum_{v,w \in V} \bar{d}(v, w)^2 \tilde{T}(v, w)\mu_0(w),$$

$$\text{s.t. } \tilde{T}\mu_0 = \mu_1,$$

assuming that a minimum exists. However, it is easy to see that computing $\tilde{T}$ results in an optimization problem whose number of variables scale as $|V|^2$, making the problem intractable for most systems of interest. Moreover, this approach requires pre-computation of all pairwise shortest path distances $\bar{d}(v, w)$, which would be computationally expensive if edge-weights are chosen to be non-uniform to incorporate possibly spatially varying perturbation costs. Hence, rather than computing an 'all-to-all' transport map $\tilde{T}(v, w)$ in one shot, one can instead use the concept of advection to compute continuous-time flow over edges $E$ to ease the computational load.

A continuous-time advection on such a graph can be described [50, 51] as,

$$\frac{d}{dt}\mu(t, v) = \sum_{e=w \to v} U(t, e)\mu(t, w) - \sum_{e=v \to w} U(t, e)\mu(t, v), \tag{12}$$

where $\mu(t, v)$ is the measure supported on the vertex $v$, and $U(t, e)$ is the flow on the edge $e$. Here we use the notation $e = v \to w$ to represent the edge $e$ directed from a vertex $v$ to $w$. Eq. (12) is the evolution equation of a continuous-time Markov chain, with the right hand side of the equation written in terms of flows over edges. This equation is graph analogue of the divergence form of advection equation in Eq. (3), and is measure preserving by construction [51]. The notion of optimal transport has been extended to such a continuous-time discrete-space setting recently [52, 53]. Following [53], one can formulate an optimal transport problem on $G$ as follows. First, define an advective inner product between two flows $U_1, U_2$ as

$$\langle U_1, U_2 \rangle_\mu = \sum_{e=v \to w} \left( \frac{\mu(v)}{\mu(w)} \cdot \frac{\mu(v) + \mu(w)}{2} \right) U_1(e) U_2(e). \tag{13}$$

Then the corresponding optimal transport distance between a set of measures $(\mu_0, \mu_1)$ supported on $V$ can be written as

$$\tilde{W}_N(\mu_0, \mu_1) = \min_{U(t,e) \geq 0} \int_0^1 \|U(t, .)\|_{\mu(t,.)} dt, \tag{14}$$

such that Eq. (12) holds for $0 \leq t \leq 1, \mu(t, v) \geq 0$, and
$$\mu(0, v) = \mu_0(v), \mu(1, v) = \mu_1(v) \quad \forall v \in V.$$

Here $\|U(t, .)\|_{\mu(t,.)} \triangleq \sqrt{\langle U, U \rangle}_\mu$. This approach is motivated by the previously discussed Benamou-Brenier approach for optimal transport on continuous spaces. Motivated by the change of variables in Section 2.2, one can define $J(t, e) \triangleq \mu(t, v)U(t, e)$ for $e = (v \to w)$. Then, Eq. (14) can be written as

$$\tilde{W}_N(\mu_0, \mu_1)^2 = \min_{U(t,e) \geq 0} \int_0^1 \sum_{e=v \to w} \left( \frac{\mu(t, v)}{\mu(t, w)} \cdot \frac{\mu(t, v) + \mu(t, w)}{2} \right) U^2(t, e) dt \tag{15}$$

$$= \min_{J(t,e) \geq 0} \int_0^1 \sum_{e=v \to w} \left( \frac{\mu(t, v) + \mu(t, w)}{\mu(t, v)\mu(t, w)} \right) \frac{J^2(t, e)}{2} dt. \tag{16}$$

Eq. (15) represents a graph version of time integral of kinetic energy of perturbation, and should be compared with cost given in Eq. (3). This results in the following advection based convex optimization problem

$$\tilde{W}_N(\mu_0, \mu_1)^2 = \min_{J(t,e) \geq 0} \int_0^1 \sum_{e=(v \to w)} \frac{J(t, e)^2}{2} \left( \frac{1}{\mu(t, v)} + \frac{1}{\mu(t, w)} \right) dt, \tag{17}$$

$$\mu(0, v) = \mu_0(v), \mu(1, v) = \mu_1(v) \quad \forall v \in V, \tag{18}$$

$$\frac{d}{dt}\mu(t, .) = D^T J(t, .) \quad 0 \leq t \leq 1, \tag{19}$$

where $J(t, e) \triangleq \mu(t, v)U(t, e)$ for $e = (v \to w)$. The matrix $D \in \mathbb{R}^{|E| \times |V|}$ is the linear flow operator computing $\mu(w) - \mu(v)$ for each $e = (v \to w) \in E$, i.e. $D^T(i, j)$ equals $+1$ if $j$th edge points into $i$th vertex, $-1$ if $j$th edge points out of $i$th vertex, and 0 if $j$th edge is not connected to $i$th vertex. Hence Eq. (19) is a rewriting of Eq. (12) in terms of $J(t,.)$. Eq. (17) is graph analogue of cost given in Eq. (5).

The convergence properties of distance $\tilde{W}_N$ have been studied in Ref. [52]. It is shown that as $N \to \infty$,

$$\tilde{W}_N(\mu_0, \mu_1) \to W_2(\mu_0, \mu_1), \tag{20}$$

where we have used the same notation for measures on continuous and discrete spaces. This results implies that for graphs obtained by discretization of continuous spaces, the measure $\tilde{W}_N$ is a meaningful quantity to compute to capture the optimal transport cost. Furthermore, for a fixed discretization size $N$, the relation between $\tilde{W}_N$ and $W_{2,N}$ has also been studied in [52]. It can be shown that,

$$\tilde{W}_N(\mu_0, \mu_1) \leq C_1 W_{2,N}(\mu_0, \mu_1), \tag{21}$$

where $C_1$ is a constant independent of $N$. This result, combined with numerical evidence in [53], show that the distance $\tilde{W}_N$ captures the graph structure of $G$, without requiring the tedious computations needed to compute $W_{2,N}$. We drop the subscript $N$ in the rest of the paper. Fig. 2 summarizes the motivation and process of arriving at $\tilde{W}_N$ as an appropriate cost for measuring phase space transport due to discrete-time perturbations on a graph.

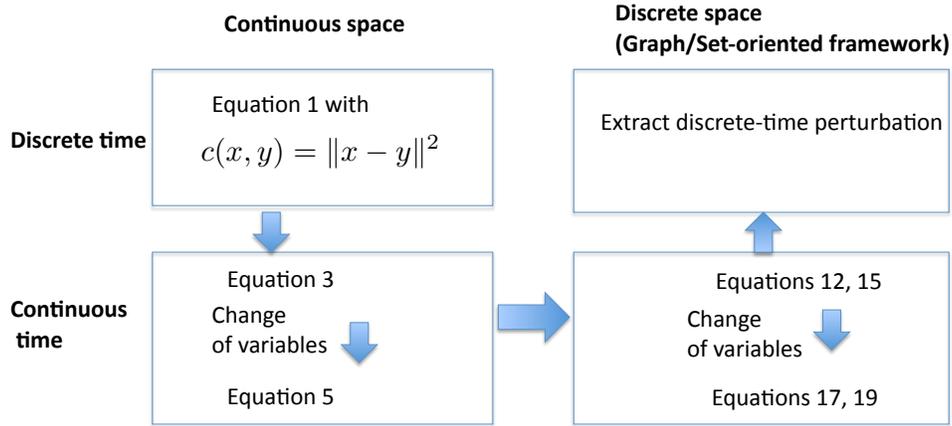

**Figure 2:** From discrete-time continuous-space formulation to discrete-time discrete-space formulation of optimal transport

Following Refs. [54, 53], we use the staggered discretization scheme for pseudo-time discretization. We define

$$\mu^j(v) \triangleq \mu(t_j, v), \tag{22}$$

$$J^j(e) \triangleq J(t_j, e), \tag{23}$$

$$\tag{24}$$

where $t_j = j/k, j \in [0, 1, 2, \ldots, k]$ is the time discretization into $k$ intervals. Here $J^j(e)$ represents the flow over edge $e : v \to w$, from vertex $v$ at $t_{j-1}$ to vertex $w$ at $t_j$.

Hence, the optimization problem given in Eqs. (17-19) can be time-discretized as,

$$\tilde{W}(\mu_0, \mu_1)^2 = \min_{J^j, \mu^j} \sum_{j=1}^{k} \sum_{\substack{e=1 \\ e=(v \to w)}}^{E} (J^j(e))^2 \left(\frac{1}{\mu^{j-1}(v)} + \frac{1}{\mu^j(w)}\right), \tag{25}$$

subject to the following constraints:

$$\frac{\mu^j - \mu^{j-1}}{\Delta t} = D^T J^j \text{ for } j \in [1, 2, \ldots, k], \tag{26}$$

$$\mu^0 = \mu_{t_0}, \mu^n = \mu_{t_f}, \tag{27}$$

where $\Delta t = \frac{1}{k}$. We note that this form of perturbation is measure preserving by construction. The cost function given by Eq. (25) is again of the form 'quadratic over linear', while the OT flow (Eq. (26)) adds linear constraints. Hence the discretized problem is convex in the problem variables $(J^j, \mu^j)$, and can be solved using many off-the-shelf convex solvers. Note that $q^j \in \mathbb{R}^{|V|}$ are vertex based, and $J^j \in \mathbb{R}^{|E|}$ are edge based quantities. The size of the optimization problem can be quantified in terms of number of vertices $|V| = N$, pseudo-time discretization $k$, and the number of edges $|E|$. The graph $G$ formed by discretization of phase space is always sparse, since a typical vertex is at most connected to $2d$ neighbors, where $d$ is the dimension of the phase space. Hence, $|E| = O(N)$, and the problem size is $O(kN)$. In practice, we see that good convergence in cost is obtained for $k \sim O(10) - O(100)$, and hence $k \ll N$ in problems of interest. Hence, this continuous (pseudo) time problem is much more tractable than solving the discrete-time problem with the cost given by $W_{2,N}$ that scales as $O(N^2)$, as discussed earlier [53]. The optimization problem is solved via CVX [55] modeling platform, an open-source software for converting convex optimization problems into usable format for various solvers. We use the SCS [56] solver, a first-order solver for large size convex optimization problems. This solver uses the Alternating Direction Method of Multipliers (ADMM) [57] to enable quick solution of very large convex optimization problems, with moderate accuracy.

Let us consider the optimal transport problem (without any system dynamics) for a configuration similar to the one considered in the original Brenier-Benamou paper [46]. The phase space is a unit 2-torus, and initial measure $\mu_0$ is taken to be uniform measure over a circular region in the center of the torus. The final measure $\mu_1$ is a translation of $\mu_0$, centered at the 'corners' of the torus. For simplicity, we use uniform measures rather than Gaussian measures considered in Ref. [46]; this does not change the solution qualitatively. The phase space is discretized into a $N = 50 \times 50$ grid. In Fig. 3, intermediate measures of the optimal transport computed via the graph based algorithm are shown. As in Ref. [46], this algorithm also results in splitting of the initial measure into four pieces, deformation, and finally, translation of each piece to the nearest corner. The convergence of optimal transport cost with $k$ is shown in Fig. 4.

### 3.3. Switching approach to optimal transport in nonlinear systems

Now we are ready to state our algorithm to solve the optimization problem given in Eqs. (6-9). As mentioned in the last section, we work with measures $\mu$ supported on $V$. Rather than solving for $\tilde{T}^i$s acting on these measures, we solve for continuous pseudo-time flows which transform the measures $\hat{\mu}^i$ to $\mu^i$. Fig. 5 describes this switching approach to solving the optimization problem. In what follows, we denote pseudo-time by $\bar{t}$, and real time by $t$.

Hence, the optimization problem written abstractly in Eqs. (6-9) can now be written as,

$$\tilde{W}_{dyn}(\mu_{t_0}, \mu_{t_f})^2 = \min_{J^i(\bar{t},e) \geq 0, q^i(\bar{t},v) \geq 0} \sum_{i=1}^{n} \int_0^1 \sum_{e=(v \to w)} \frac{J^i(\bar{t},e)^2}{2} \left( \frac{1}{q^i(\bar{t},v)} + \frac{1}{q^i(\bar{t},w)} \right) d\bar{t}, \tag{28}$$

$$\hat{\mu}^i = \mu^{i-1}\hat{P}, \tag{29}$$

$$q^i(0,v) = \hat{\mu}^i(v), q^i(1,v) = \mu^i(v) \quad \forall v \in V, \tag{30}$$

$$\frac{d}{d\bar{t}}q^i(\bar{t},.) = D^T J^i(\bar{t},.) \ 0 \leq \bar{t} \leq 1, \tag{31}$$

$$\mu^0 = \mu_{t_0}, \mu^n = \mu_{t_f}. \tag{32}$$

Here $q^i(\bar{t},.)$ is the pseudo-time varying measure, and $J^i(\bar{t},e) \triangleq q^i(\bar{t},v)U^i(\bar{t},e)$ is the flow for the $i$th instance of optimal transport. To numerically solve the optimization problem, we time-discretize the advection equation, i.e. Eq. (31), and corresponding flow $J(.,\bar{t})$, w.r.t pseudo-time variable $\bar{t}$ for each OT step.

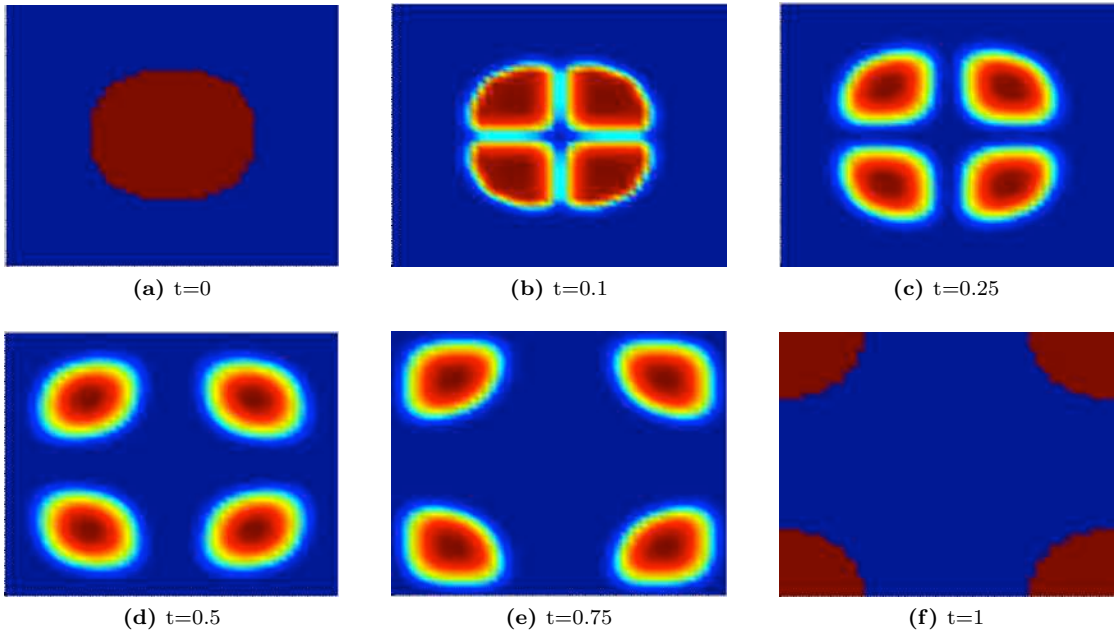

**Figure 3:** Graph-based optimal transport of uniform measure with circular support, from the center to the corner of the torus.

The staggered scheme for pseudo-time discretization discussed in Section 3.2 is now modified to incorporate evolution under system dynamics between the OT steps. We define

$$q^{i,j}(v) \triangleq q^i(\bar{t}_j, v), \tag{33}$$

$$J^{i,j}(e) \triangleq J^i(\bar{t}_j, e), \tag{34}$$

$$\tag{35}$$

where $\bar{t}_j = j/k, j \in [0, 1, 2, \ldots, k]$ is the discretization into $k$ intervals. Here $J^{i,j}(e)$ represents the flow over edge $e : v \to w$, from vertex $v$ at $\bar{t}_{j-1}$ to vertex $w$ at $\bar{t}_j$. The resulting optimization problem can be written as,

$$\tilde{W}_{dyn}(\mu_{t_0}, \mu_{t_f})^2 = \min_{J^{i,j}, q^{i,j}} \sum_{i=1}^{n} \sum_{j=1}^{k} \sum_{\substack{e=1 \\ e=(v \to w)}}^{E} (J^{i,j}(e))^2 \left(\frac{1}{q^{i,j-1}(v)} + \frac{1}{q^{i,j}(w)}\right), \tag{36}$$

subject to the following constraints:

$$\hat{\mu}^i = \mu^{i-1} \hat{P}, \tag{37}$$

$$q^{i,0}(v) = \hat{\mu}^i(v), q^{i,k}(v) = \mu^i(v) \quad \forall v \in V, \tag{38}$$

$$\frac{q^{i,j} - q^{i,j-1}}{\Delta \bar{t}} = D^T J^{i,j} \text{ for } j \in [1, 2, \ldots, k], \tag{39}$$

$$\mu^0 = \mu_{t_0}, \mu^n = \mu_{t_f}, \tag{40}$$

where $\Delta \bar{t} = \frac{1}{k}$. The cost function given by Eq. (36) is again of the form 'quadratic over linear', while the evolution due to dynamics (Eq. (37)), and OT flow (Eq. (39)) are both linear constraints. The problem variables for the convex optimization problem Eqs. (36-40) are vertex based quantities $q^{i,j} \in \mathbb{R}^{|V|}$, and edge based quantities $J^{i,j} \in \mathbb{R}^{|E|}$. The variables in this optimization problem are $O(kn(N + |E|)) = O(knN)$.

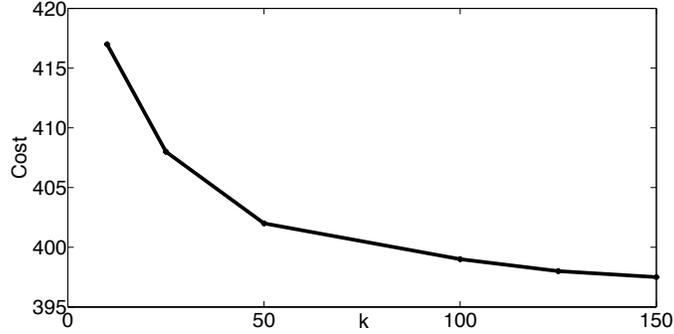

**Figure 4:** Convergence of optimal transport cost Eq. (25) for the Brenier-Benamou system with time discretization $k$.

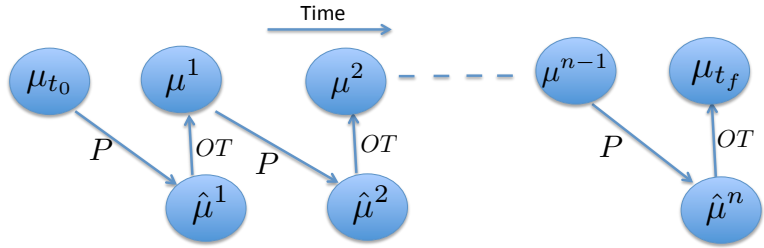

**Figure 5:** The switching approach to optimal transport for nonlinear systems. The transfer operator $P$ pushes forward measures $\mu^{i-1}$ to $\hat{\mu}^i$ over one time period. The optimal transport based perturbation step, labeled 'OT', and carried out in pseudo-time, pushes forward $\hat{\mu}^i$ to $\mu^i$.

## 4. Examples

### 4.1. Revisiting the Brenier-Benamou example: With linear dynamics

We first demonstrate our algorithm on a simple modification of the transport problem on torus, which was considered in Section 3. We consider the following discrete-time linear dynamics on the phase space $X = \mathbb{T}^2$.

$$x_{t+1} = x_t + 0.2, \tag{41}$$
$$y_{t+1} = y_t. \tag{42}$$

This problem is solved for $n = 5$ time steps, and $k = 20$ pseudo-time discretization for representing each perturbation with the OT step. Adding the dynamics breaks the symmetry of the system, and that of the resulting transport. In Fig. 6, the resulting transport is shown. The dynamics and final time are chosen such that the transport in $x$ direction can be solely accomplished by the just 'going with flow', i.e. without any perturbations. It is clear from Fig. 6 that the resulting transport indeed makes use of this fact, and the perturbations by OT only occur in the $y$ direction.

### 4.2. Optimal transport in time-periodic double-gyre system

Now we consider a measure transport problem for the time-periodic double-gyre system [58]. This system is well-studied in the dynamical systems literature, and it has been extensively used as a proving ground for several new computational tools related to transport and mixing [58, 48, 59, 39].

The system equations are as follows,

$$\dot{x} = -\pi A \sin(\pi f(x,t)) \cos(\pi y), \tag{43}$$
$$\dot{y} = \pi A \cos(\pi f(x,t)) \sin(\pi y) \frac{df(x,t)}{dx}, \tag{44}$$

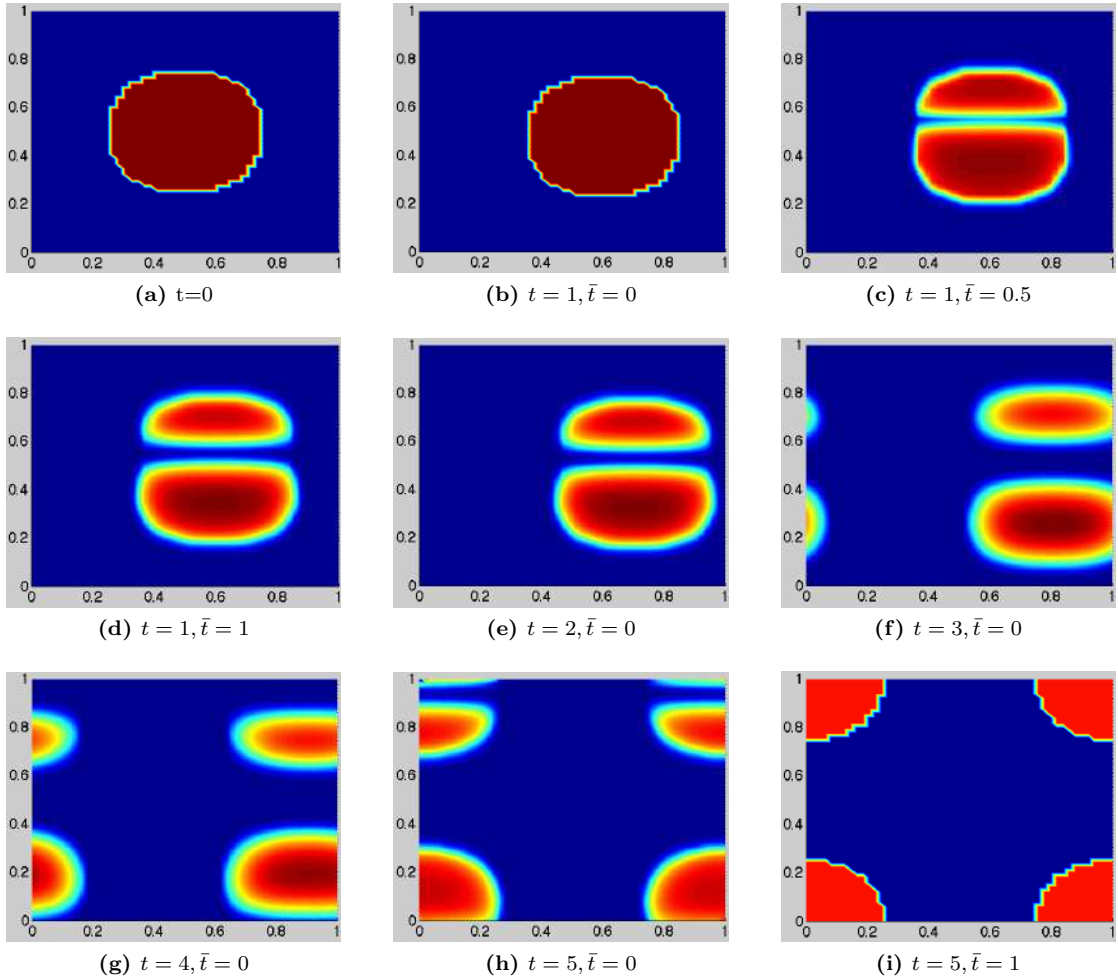

**Figure 6:** Graph-based optimal transport with dynamics: Transport of uniform measure with circular support, from the center to the corner of the torus with dynamics given by Eqs. (41-42).

where $f(x,t) = \beta \sin(\omega t) x^2 + (1 - 2\beta \sin(\omega t))x$ is the time-periodic forcing in the system. The phase space is $X = [0\ 2] \times [0\ 1]$. The dynamics of this system can be understood by a combination of geometric and statistical tools. For the trivial case of $\beta = 0$ (i.e. no time-dependent forcing), the phase space is divided into two invariant sets, i.e., the left and right halves of the rectangular phase space ('gyres'), by a heteroclinic connection between fixed points $x_1 = (1,1)$ and $x_2 = (0,1)$. For non-zero $\beta$, it is instructive to study the dynamics of Poincare map $F$ of the system, obtained by integrating the dynamics over one time period $\tau$ of $f$. The heteroclinic connection is broken in this case, and results in a heteroclinic tangle. This heteroclinic tangle leads to transport between left and right gyres via lobe-dynamics.

In this study, we use $A = 0.25, \beta = 0.25, \omega = 2\pi$, which results in time-period $\tau = 1$. The theory of lobe dynamics describes the qualitative and quantitative aspects of inter-gyre transport [58]. In figure 7, the unstable manifold of $x_1 \approx (0.919, 1)$, $U_{x_1}$, and the stable manifold of $x_2 \approx (1.081, 0)$, $S_{x_2}$ are shown in green and white respectively. The lobe labeled 'A', its pre-image $F^{-1}(A)$ and image $F(A)$ are also shown. Consider the segment $L = U_{x_1}(x_1 \to F(P_1)) \cup S_{x_2}(F(P_1) \to x_2)$. Then $L$ divides phase space $X$ into two regions. The points that get mapped from left to right region in one iteration of $F$ are precisely those in set $A$. Hence, the amount of mass transport from left to right side of $L$ is $\bar{m}(A)$.

For further insight into transport in this system, set-oriented transfer operator methods have been utilized

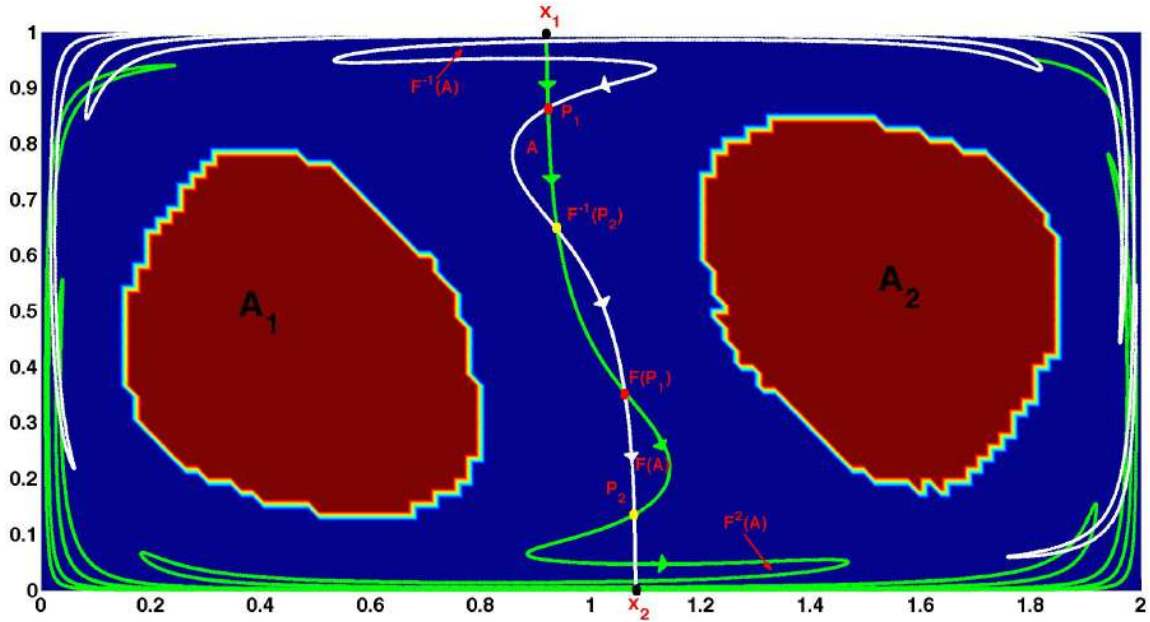

**Figure 7:** Invariant manifolds and lobe-dynamics in the double-gyre system.

previously in the literature [48]. A set $S$ is called almost-invariant if

$$\frac{\bar{m}(F^{-1}(S) \cap S)}{\bar{m}(S)} \approx 1.$$

Invariant and almost-invariant sets in this system can be identified by the sign structure of the second eigenvector of the reversibilized transfer operator

$$R = \frac{P + \tilde{P}}{2},$$

where $\tilde{P}$ is the transfer operator corresponding to reverse-time dynamics. In Fig. 7, two almost-invariant sets, $A_1$ and $A_2$ are also shown.

We take $\mu_0$ to be the uniform measure supported on $A_1$, and $\mu_{t_f}$ to be the uniform measure supported on $A_2$. Both measures are normalized to sum to unity. We solve the the optimal transport problem for $N = 100 \times 50$ & $k = 10$ for different time horizons. We use finer grid-sizes $N = 120 \times 60$ and $N = 150 \times 75$ to verify that our results are nearly independent of $N$. Recall that $t_f = n\tau$, where $\tau$ is the period of the flow.

In Fig. 8, the solution with $n = 1$, i.e., $t_f = \tau = 1$ is shown. Recall that $0 \leq t \leq t_f$ is time, and $0 \leq \bar{t} \leq 1$ is used to denote pseudo-time variable used in computation of the discrete-time perturbations. We use $\bar{t} = 0$ to denote beginning of each perturbation step, and $\bar{t} = 1$ for end of the perturbation step. Since for $n = 1$ there is only one time-step of dynamics and one discrete-time perturbation, most of the transport happens in the perturbation step. The result of sole time-step of the dynamics is shown in Fig. 8(b), which is similar to Fig. 8(a) since the initial measure is supported on AIS $A_1$. In Figures 8(b)-(f), we show the pseudo-time evolution of the sole discrete-time perturbation which begins at $t = 1, \bar{t} = 0$ and ends at $t = 1, \bar{t} = 1$.

We show the solution of the optimization problem with larger time-horizon $t_f = 15$ in Fig. 9. In this case, a large fraction of the measure is pushed to the right via lobe-dynamics, leading to efficient transport. The role of perturbations for large time-horizons is to push the measure into these lobes (and their pre-images). This concentration will ensure that lobe-dynamics pushes measure from left to right gyre. For instance, at the end of third perturbation, shown in Fig. 9(b), mass has been pushed close to the sets $F^{-1}(A)$ and $A$.

The next application of the dynamics lead to this mass being pushed to the neighborhood of $A$ and $F(A)$, respectively, as shown in Fig. 9(c). This pattern is repeated over next several time-steps, eventually leading to 'emptying' out of $A_1$, and collection of mass on $A_2$. Note that since we are solving a fixed-final time problem, the intermediate measures can be very different from the final target measure. The optimization problem allows for arbitrary mixing of intermediate measures, and optimally concentrates it on regions which will be transported to the right side over the next iteration. The process of concentration itself is carried out within perturbation pseudo-time steps, as can been seen by comparing (say) Fig. 9(c) (beginning of 4th perturbation) with Fig. 9(d) (end of 4th perturbation).

In Fig. 10, we show the cost of each discrete-time perturbation over the time horizon of the problem. The use of efficient global transport given by lobe dynamics results a drastic decrease in value of the optimal transport cost $\tilde{W}$, as shown in Fig. 11.

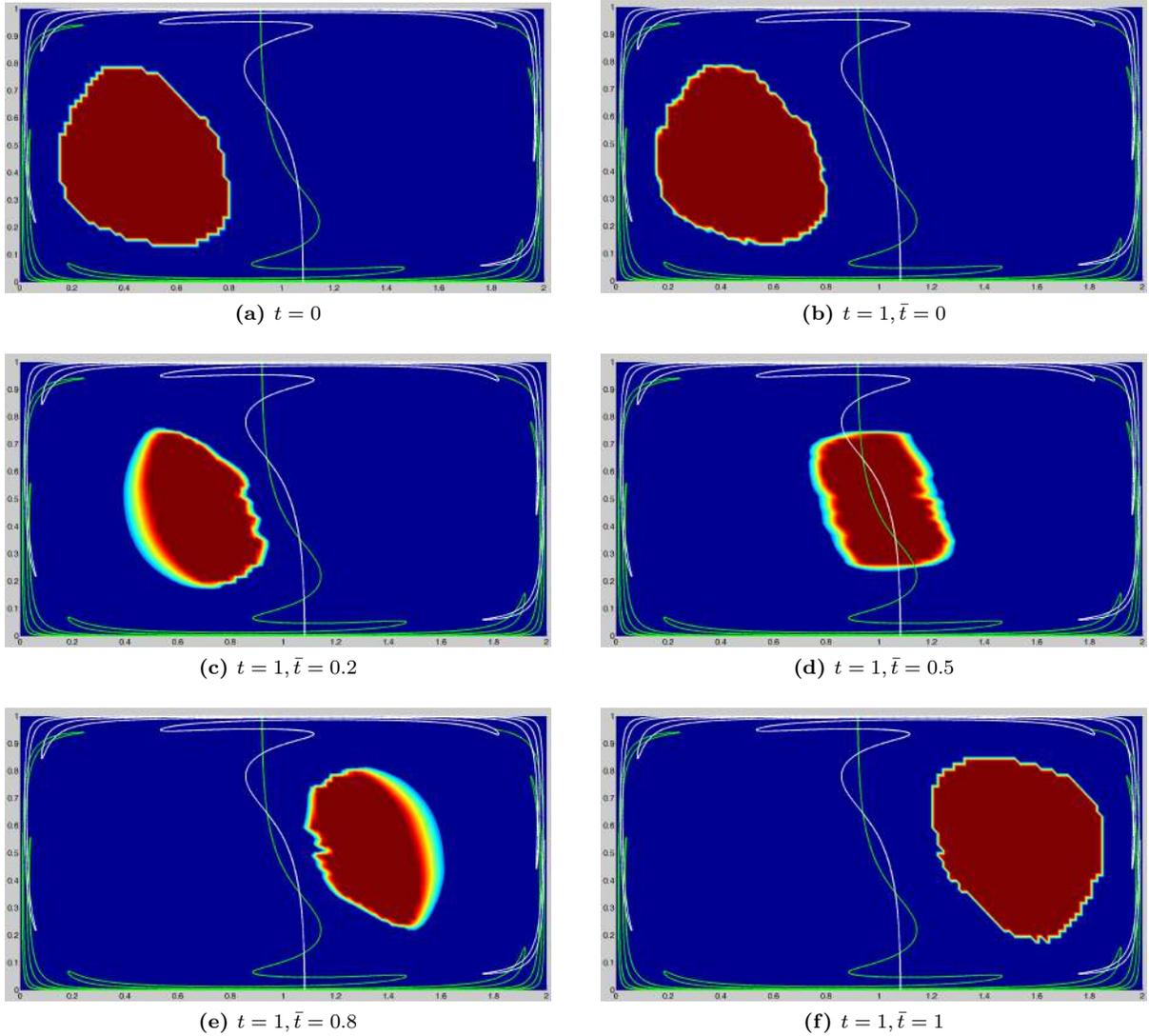

**Figure 8:** Graph-based optimal transport in time-periodic double-gyre: Transport of a measure from the left AIS to the right AIS, with $t_f = \tau = 1$. Notice since there is only one perturbation (OT) step, all of the initial measure is pushed across the invariant manifolds to the final measure by the pseudo-time OT flow.

### 4.2.1. Localization of Perturbations

Since there exist pathways for global transport due to lobe dynamics, we expect the perturbations to increasingly localize as time-horizon $t_f$ is increased. In other words, we expect that for large time horizons, globally optimal perturbations would not involve transport which can accomplished for 'free'. To quantify this behavior, we compute the map $\tilde{T}^i : \mathbb{R}^{|V|} \to \mathbb{R}^{|V|}$ corresponding to $i$th perturbation in the optimization problem given by Eqs. (36-40). We first reconstruct the advection map over each pseudo-time step in $i$th perturbation, $\tilde{T}^{i,j} : \mathbb{R}^{|V|} \to \mathbb{R}^{|V|}$, using an infinitesimal-generator approximation [60], as follows.

$$\tilde{T}^{i,j}(v,w) = \begin{cases} \sum_{e:w \to v} U^{i,j}(e) & v \neq w \\ 1 - \sum_{e:v \to \hat{v}} U^{i,j}(e) & \text{otherwise} \end{cases} \quad (45)$$

Then, $\tilde{T}^i$ is obtained as the composition map of $k$ pseudo-time advection maps, i.e $\tilde{T}^i = \tilde{T}^{i,k}\tilde{T}^{i,k-1}\ldots\tilde{T}^{i,1}$. We define a distance $D_l$ on the set of all $\tilde{T}^i, i = 1, 2, n$ as follows.

$$D_l = \frac{1}{\delta x} \max_{i \in [1,2,\ldots,n]} \max_{v \in V} \max_{\{w | \tilde{T}^i(v,w) > 0\}} d_\infty(v,w), \quad (46)$$

where $d_\infty(v,w)$ is the greater of the distances in $x$ and $y$ directions between centers of boxes $v$ and $w$. Hence $D_l$ is simply the greatest distance in $x$ or $y$ directions that any mass is moved by any of the $n$ perturbations, divided by $\delta x$. Here $\delta x$ is the length of $B_i$. In Fig. 12, we plot the measure $D_l$ for various values of $t_f$. It is clear from this plot that along with decreasing transport cost, the perturbations tend to localize with increasing $t_f$. However, there seems to be a minimum value of $D_l \approx 8$, even for large $t_f$, indicating a lower limit on the 'radius' of perturbation for reaching the desired final measure.

### 4.3. Optimal mixing enhancement in time-periodic Double-Gyre system

Next, we apply our algorithm to the problem of optimal enhancement of finite-time mixing. We aim to obtain optimal perturbations at discrete times that lead to complete mixing in given finite time. We choose $\mu_{t_0}$ to be the uniform measure supported on $A_1$, and $\mu_{t_f}$ to be the uniform measure over phase space $X$. Both measures are normalized to unity as before.

Note again that this problem setting is different than those considered in previous works on mixing measures and optimal mixing [21, 20, 39], since the aim is to minimize a cost associated with phase space transport due to applied perturbations, rather than maximizing 'mixedness' at final time. Furthermore, since we are using a set-oriented framework, there is an inherent minimal scale present in the problem, i.e. the size of each box $B_i$ or vertex $v$. Inhomogeneity in the measure at length scales below this size are ignored, and hence, complete mixing in this context implies that inhomogeneities for all scales above the size of smallest box in the partition have been removed. We also note that using set-oriented approach for propagating measures is known to cause artificial diffusion.

The solution sequence for $n = 1$ is shown in Fig. 13. In this short time-horizon mixing problem, only some of the mass is transported directly across the invariant manifolds to the right side, in contrast with the case of $n = 1$ transport problem considered earlier. The rest of the mass is mixed on the left side by a smearing type action, and the mass transported to right is similarly mixed. For the long time-horizon case shown in Fig. 14, the situation is analogous. Some of the perturbation effort is spent in smearing the measure on each side, while rest of time is used to push it into the lobes (and their pre-images) to enable global transport of 'half' of the measure to the the right side. Both these actions are carried out simulataneously. For instance, the third perturbation places some mass near $F^{-1}(A)$ and $A$, while mixing other mass all across the left side, shown in Fig. 14(c). The action of $P$ moves mass to $A$ and $F(A)$ respectively, as shown in Fig. 14(d). This process is again repeated over the rest of the time-horizon. Fig. 14(g) shows the final lobe transport from left to right side, while the rest of the measure is almost fully mixed. In Fig. 15, we show the cost of each discrete-time perturbation over the time horizon of the problem. The use of efficient global transport given by lobe dynamics again results a large decrease in value of the optimal mixing cost $\tilde{W}$, as shown in Fig. 16.

In Fig. 17, we plot the measure $D_l$ for various values of $t_f$. It can be seen that along with decreasing mixing cost, the perturbations tend to localize with increasing $t_f$. Compared to the transport case, the optimal mixing solution shows higher localization with increasing $t_f$. Hence, the localized perturbations play a dual role. On one hand, they enable the use of lobe dynamics for global transport by concentrating measure in pre-images of lobes (the measure concentration can be interpreted as a de-mixing process). On the other hand, they increase local mixing which leads to homogenization at a faster rate than that achieved by natural dynamics alone. The localization information obtained from computing $D_l$ can be valuable in designing efficient short-time mixing protocols, and formulating an optimization problem with localization constraints on the perturbations.

## 5. Conclusions and Extensions

We have presented a framework for computing provably globally optimal perturbations for transporting a given initial phase space measure to an arbitrary final measure in finite time for nonlinear dynamical systems. We model the discrete-time perturbations as maps which result from Monge-Kantorovich optimal transport on graphs. Our work represents a first step towards combining the set-oriented transfer operator methods with optimal transport theory, via use of graph-based optimal transport concepts. Our results show that the resulting transport exploits efficient global transport available via mechanisms such as lobe-dynamics, as the time-horizon of the problem is increased.

Several extensions of this work are desirable. The extension to dynamical systems with arbitrary time dependence is straightforward. The interesting case to study in these systems is optimal transport between coherent sets. By exploiting efficient phase space discretization techniques, such as those employed in GAIO [9], one can hope to improve the efficiency of the resulting optimal transport algorithms, and apply the framework to higher dimensional dynamical systems. Graph pruning algorithms can be employed to remove edges which are not likely to be used during the perturbation step [61].

While we give numerical evidence that the perturbation cost decreases with increasing time-horizon, and that the perturbation becomes increasingly localised, it is desirable to obtain more quantitative estimates of such behavior, possibly via a rigorous convergence analysis of the problem in the long-time limit. We have not explictly considered control constraints in the current work. By using edge weights on edges $E$ which reflect some control cost, this framework can potentially be extended to obtain perturbations which can be implemented via available control mechanisms, and result in optimal control cost. This is a challenging problem since it involves deriving conditions on controllability or reachability of measures in discrete-time setting. One way to make progress would be relax the final condition to only approximately achieve the final distribution. For chaotic systems, this relaxed version should result in a well-posed problem under mild assumptions. Connections with work in the closely related area of occupation measures [31] and Lyapunov measures [28] also need to be explored. Comparing optimal mixing analysis given in this paper with prior work which has exploited a different optimal transport measure to obtain bounds on finite time mixing [22] is another topic of possible interest.

Recent methods in obtaining Lagrangian coherent structures and coherent sets in finite-time non-autonomous systems have used variational formulations of transport under nonlinear dynamics [14, 62]. It would be fruitful to develop connections of these objects with optimal mass transportation theory, since there already exist such connections in the autonomous Hamiltonian dynamics case [43].

Finally, it is of great interest to obtain a set-oriented approach to continuous-time optimal transport formulation of nonlinear systems with control vector fields, i.e. obtaining control laws that implement possibly non-holonomic optimal transport of measures in nonlinear dynamical systems. This topic is considered in a forthcoming work.

## 6. Acknowledgments

We thank the anonymous reviewers for their careful reading of the manuscript, and helpful suggestions that have led to significant improvement in the quality of this work.

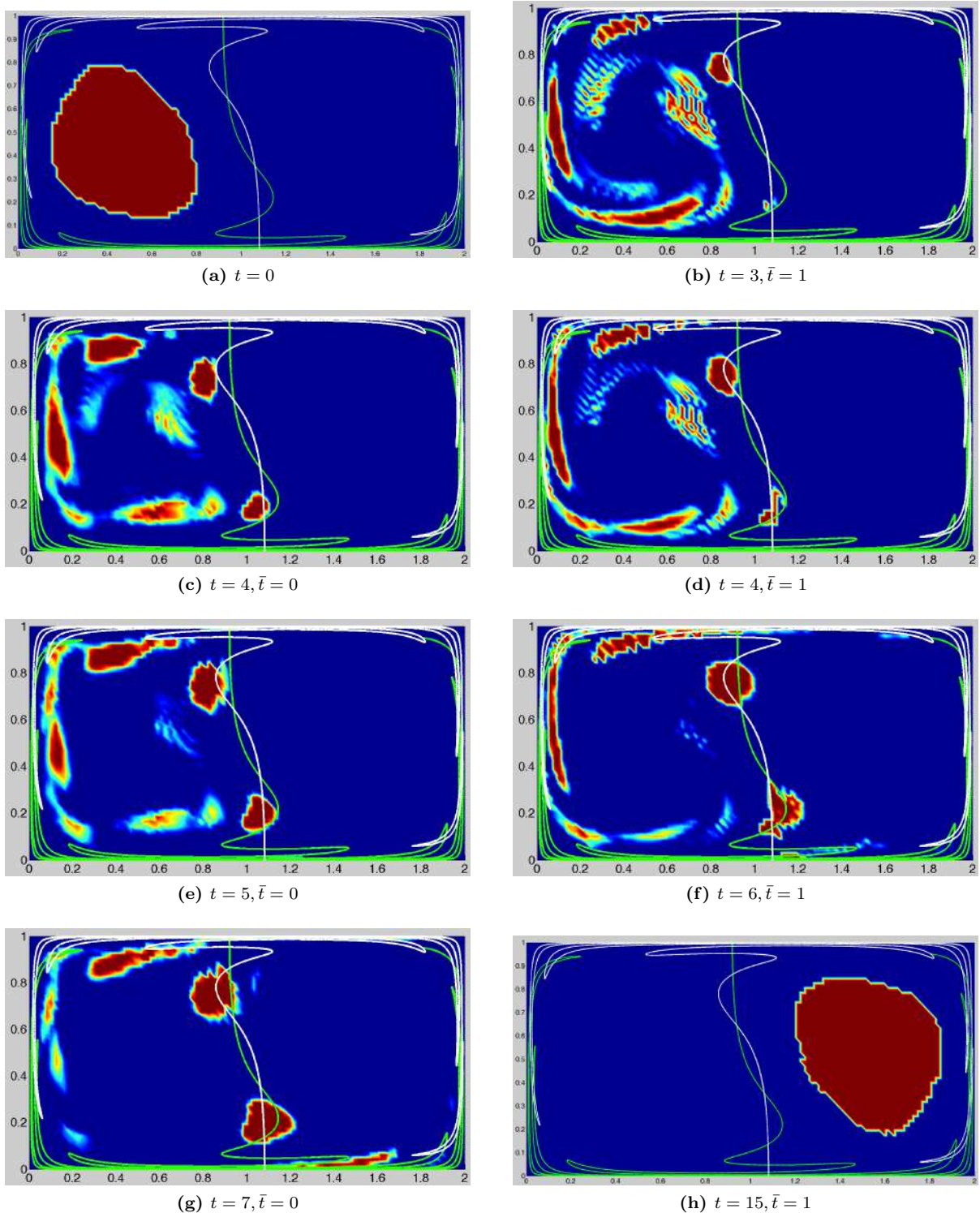

**Figure 9:** Graph-based optimal transport in time-periodic double-gyre: Transport of a measure from the left AIS to the right AIS, with $t_f = 15\tau = 15$. (a) The initial measure. (b). Perturbation pushes measure towards $F^{-1}(A)$ and $A$. See Fig. 7 for notation. (c) An iteration of $F$ maps measure from neighborhood of $A$ to neighborhood of $F(A)$. (d) Perturbation pushes measure into $F^{-1}(A)$, $A$ and $F(A)$. (e) An iteration of $F$ maps measure at $F^{-1}A$, $A$ into $A, F(A)$ respectively. (f) Perturbation places measure at $F^{-1}(A), A, F(A)$. (g) Next iteration of $F$ moves the same measure into $A, F(A), F^2(A)$. The sequence shown in (f)-(g) for $t = 6$ repeats itself over several time-steps. (h) The final measure is supported on $A_2$. **Animation:Video file NX100_k10_n15_transport.mpg**

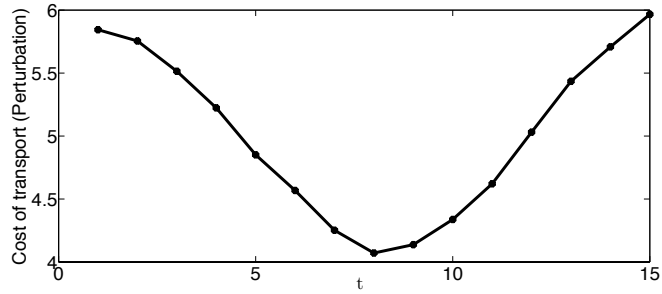

**Figure 10:** The cost of discrete-time perturbations $T^i$, $i = 1$ to 15, for the case shown in Fig. 9.

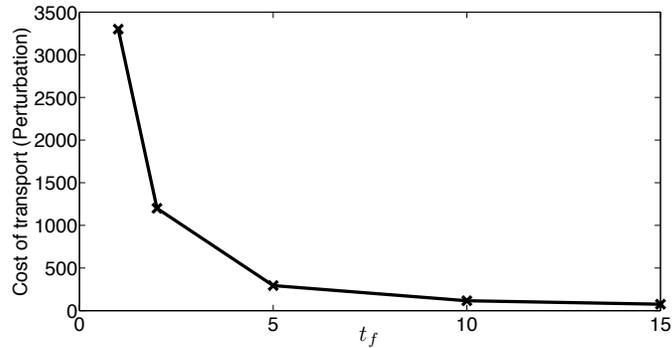

**Figure 11:** The optimal transport cost (given by Eq. 36) for the double-gyre system for transport between the two measures with supports $A_1$ and $A_2$ respectively, for various time-horizons $t_f$. For small time-horizons, perturbations push bulk of the measure directly into the right side. As the time-horizon increases, a larger fraction of the mass is pushed to the right via lobe-dynamics, leading to efficient transport. The role of perturbations for large time horizons is to push the measure into these lobes (and their pre-images).

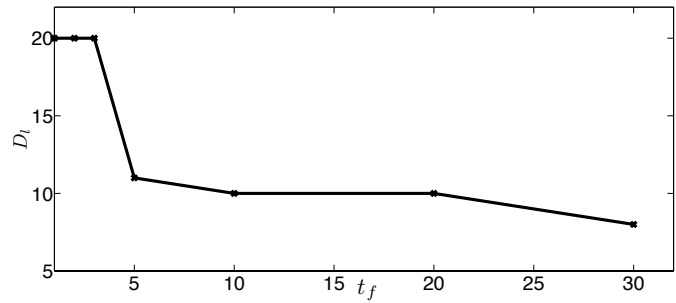

**Figure 12:** The measure $D_l$, defined in Eq. (46), as a function of time-horizon $t_f$. Decay in the the value of $D_l$ implies that perturbations get increasingly localized as time-horizon is increased

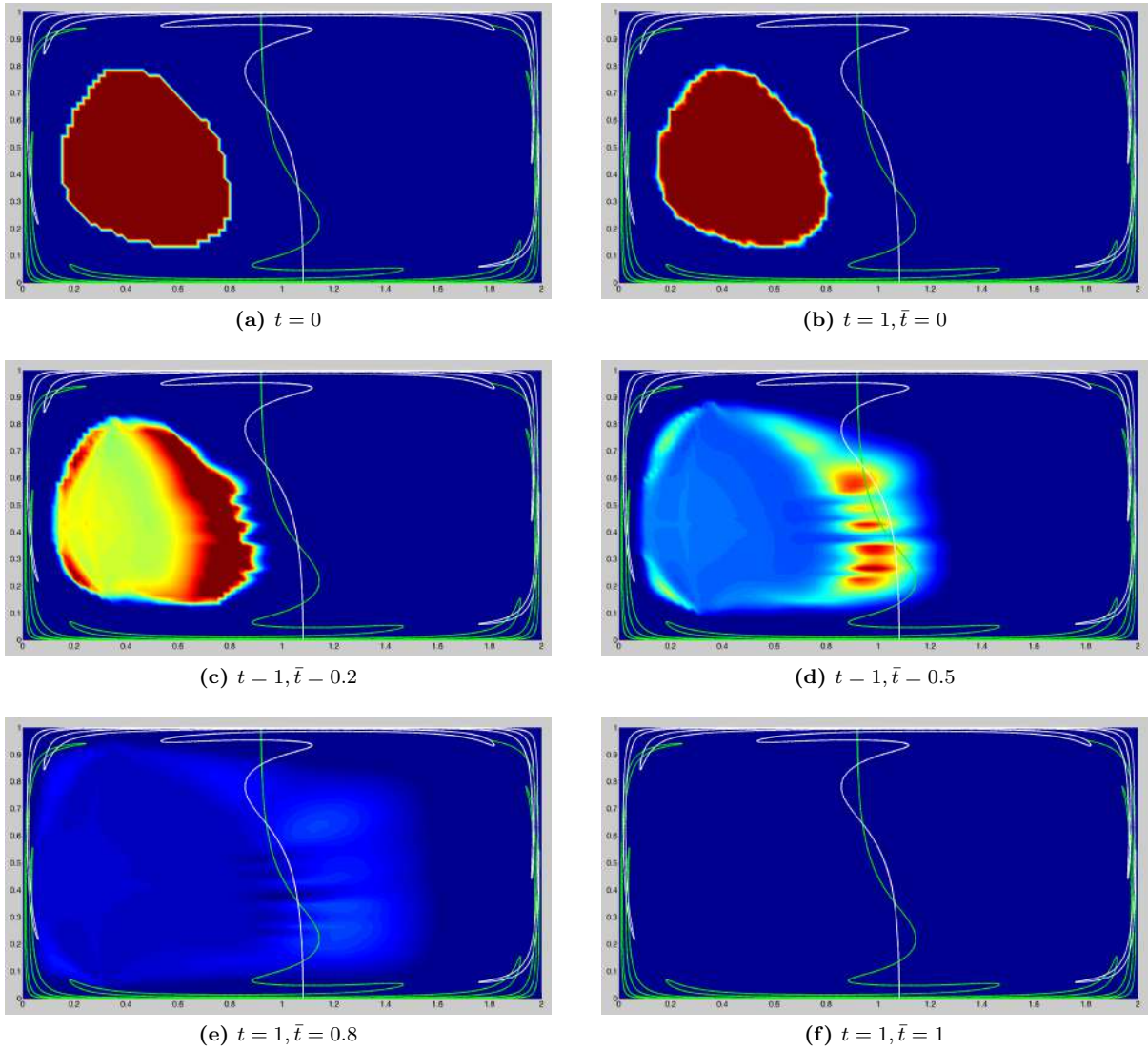

**Figure 13:** Graph-based optimal mixing in time-periodic double-gyre with $t_f = \tau = 1$. Notice since there is only one perturbation (OT) step, all of the initial measure is pushed across the invariant manifolds to the final measure by the pseudo-time OT flow).

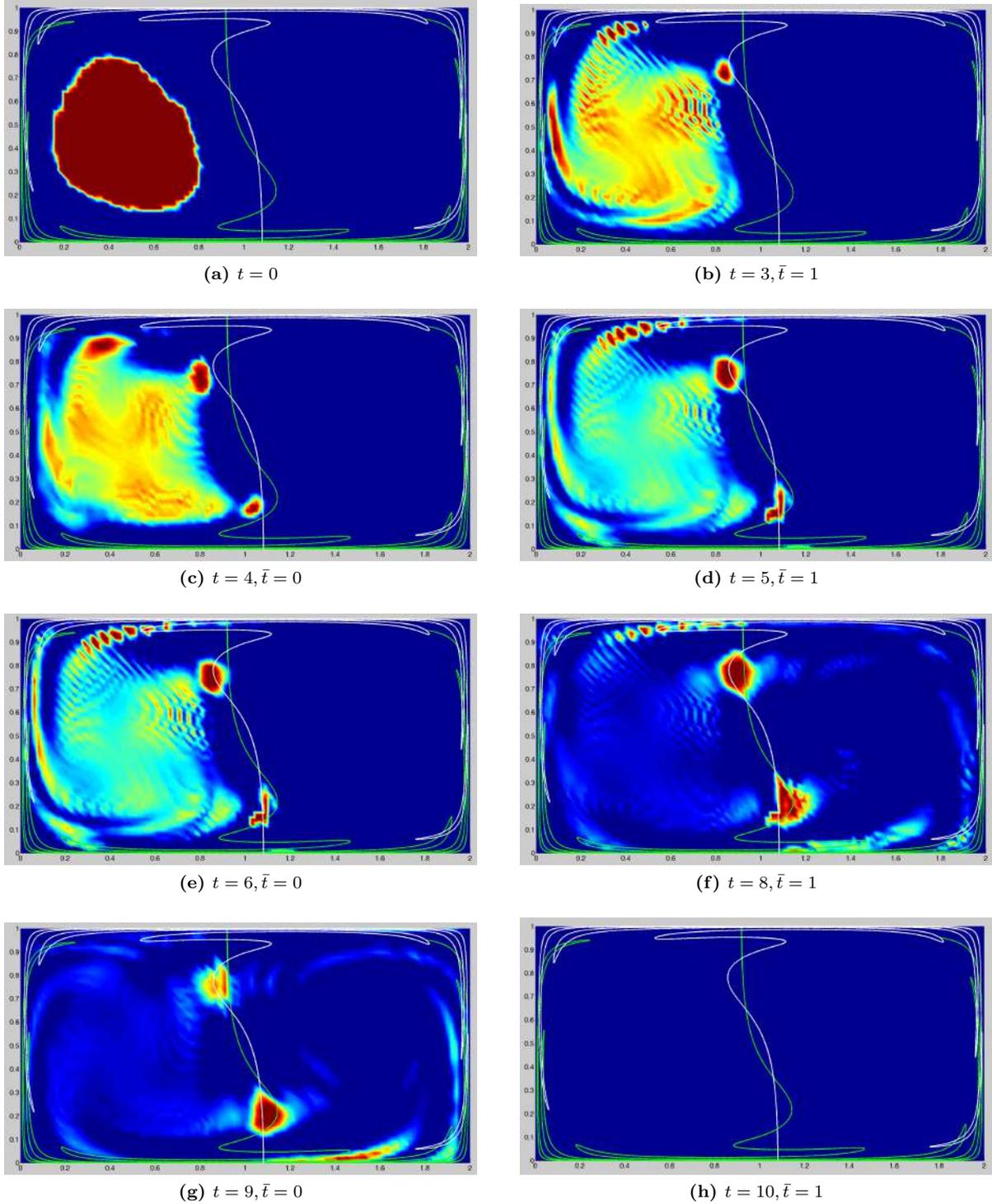

**Figure 14:** Graph-based optimal mixing in time-periodic double-gyre with $t_f = 10\tau = 10$. (a) The initial measure. (b). Perturbation pushes measure towards $F^{-1}(A)$ and $A$. See Fig. 7 for notation. (c) An iteration of $F$ maps measure near $F^{-1}(A)$ and $A$ to near $A$ and $F(A)$, respectively. (d) Perturbation pushes measure into $F^{-1}(A)$, $A$ and $F(A)$. (e) An iteration of $F$ maps measure at $F^{-1}A, A$ into $A, F(A)$ respectively. (f) Perturbation places measure at $F^{-1}(A), A, F(A)$. (g) Next iteration of $F$ moves the same measure into $A, F(A), F^2(A)$. The sequence shown in (f)-(g) for $t=6$ repeats itself over several time-steps. (h) The final measure is uniform. **Animation : Video file NX100_k10_n10_mixing.mpg**

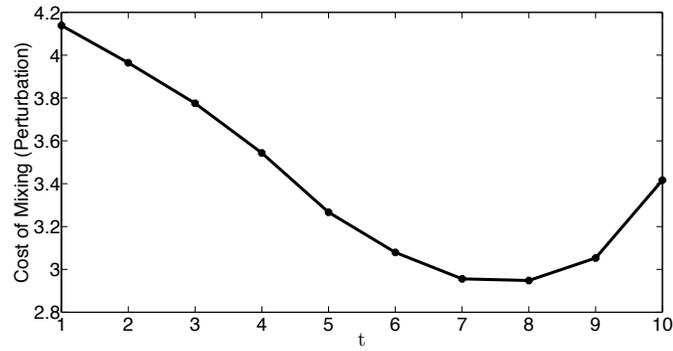

**Figure 15:** The cost of discrete-time perturbations $T^i$, $i = 1$ to $10$, for the case shown in Fig. 14.

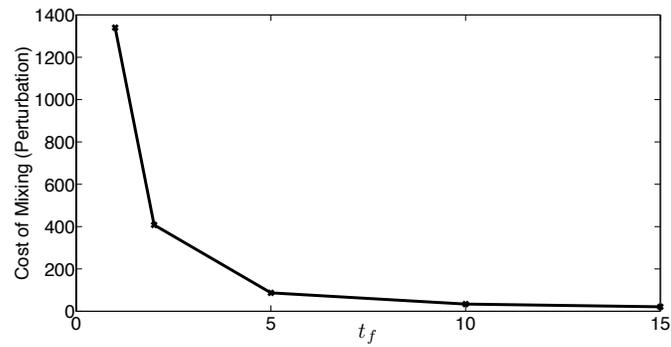

**Figure 16:** The optimal mixing cost (given by Eq. 36) for the double-gyre system for mixing the measure supported on $A_1$, for various time-horizons $t_f$. For small time-horizons, perturbations push half of the measure directly into the right side, before smearing it to make final measure uniform. As the time-horizon increases, a larger fraction of the mass is pushed to the right via lobe-dynamics, leading to drastically lower mixing cost. The role of perturbations for large time horizons is to push the measure into these lobes (and their pre-images), and later smear it uniformly in the domain.

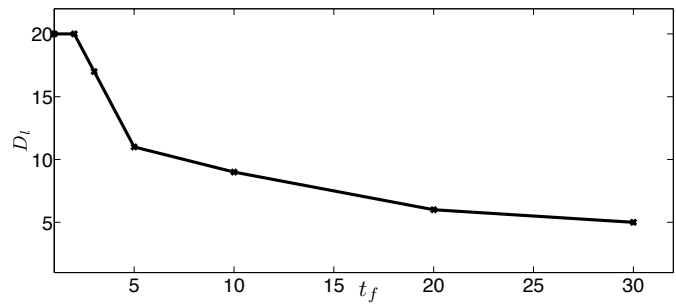

**Figure 17:** The measure $D_l$, defined in Eq. (46), as a function of time-horizon $t_f$. Decay in the the value of $D_l$ implies that perturbations get increasingly localized as time-horizon is increased.